\documentclass[a4paper,10pt]{article}
\usepackage[english]{babel}
\usepackage[latin1]{inputenc}
\usepackage{amsmath,amssymb,amsthm, fancyhdr,latexsym}
\usepackage{eurosym}
 \usepackage{avant}
 \usepackage{algorithm,algpseudocode}
\usepackage{graphicx}
\usepackage{subfigure}
\usepackage{color}
\usepackage{tikz}
\usetikzlibrary{arrows,shapes,matrix,decorations.pathmorphing}
\usetikzlibrary{backgrounds}
\newtheorem{theorem}{Theorem}[section]
\newtheorem{proposition}[theorem]{Proposition}

\newtheorem{example}[theorem]{Example}

\newtheorem{corollario}[theorem]{Corollary}

\newtheorem{definition}[theorem]{Definition}

\newtheorem{lemma}[theorem]{Lemma}
\newtheorem{lemmadef}[theorem]{Lemma $/$ Definition}
\newtheorem{oss}[theorem]{Remark}

  \newenvironment{dimo}{ \noindent {\bf Proof:} }{$\diamond$\vspace{0.5cm}} 
     
     \newenvironment{prob}{ \noindent {\bf Problem:} }{$\circledast$ \vspace{0.5cm}}        
             
\makeindex

\newcommand{\NN}{\mathbb{N}}

\newcommand{\CC}{\mathbb{C}}

\newcommand{\cemu}{Cerlienco-Mureddu \;}
\newcommand{\st}{ such that\;}

\newcommand{\gb}{ Groebner basis \;}
\newcommand{\gbx}{ Groebner bases \;}

\newcommand{\pe}{ potential expansion \;}
\title{A simplified version of the  `Axis of Evil Theorem' for distinct points.}
\author{{Michela Ceria}\\{Universit\`{a} degli Studi di Torino.}\\ {\tt michela.ceria@unito.it}}
\date{}
\begin{document}
\fussy
\maketitle
\begin{abstract}
Given a finite set $\mathbf{X}$ of distinct points, Marinari-Mora's `Axis of Evil Theorem' states that a combinatorial algorithm and interpolation enable to find a `linear' factorization for a lexicographical minimal \gb $\mathcal{G}(I(\mathbf{X}))$ of the zerodimensional radical ideal $I(\mathbf{X})$.\\
In this work we provide such algorithm, showing that it ends in a finite number of steps and that it actually provides the correct result.\\
The `Axis of Evil' algorithm takes as input the monomial basis of the initial ideal $T(I(\mathbf{X}))$ but its starting point is the (finite) Groebner escalier $N$ (obtained via \cemu correspondence) so we will also define the `potential expansion' 's algorithm, a combinatorical algorithm which computes the minimal basis from a finite Groebner escalier.
\end{abstract}
\section{Introduction.}\label{Intro}
Marinari-Mora in  \cite{IRAN}, \cite{CUBA}, \cite{ISSAC} gave a deep description of the structure of a  zero-dimensional ideal $I$ described by giving its Macaulay basis $\mathcal{B}(I)$ (\cite{SPES2}); in particular they enhanced the description of the  Grobner basis of an ideal in $K[X,Y]$ given by Lazard in \cite{La} proving that in a restricted case which includes the radical one, for each monomial  $\tau :=  X_1^{d_1}\cdots X_n^{d_n}$ belonging to the minimal basis $G(I)$ of the initial ideal of $I$, it is possible to produce  linear factors $\gamma_{m\delta\tau}:=X_i-f(X_1,\ldots,X_{i-1}), 1\leq m \leq  n, 1\leq\delta\leq d_m$ such that the polynomials
 $f_\tau := \prod_{m=1}^n \prod_{\delta=1}^{d_m}\gamma_{m\delta\tau}$ form a
 minimal lexicographical Groebner basis of $I$;  each such factors were obtained by producing an appropriate decomposition of the given Macaulay basis $\mathcal{B}(I)=\bigsqcup_{m=1}^n \bigsqcup_{\delta=1}^{d_m} {\sf S}_{m\delta}(\tau)$ 
and interpolating over the  monomial set obtained applying Cerlineco-Mureddu Algorithm over the set of
 functionals ${\sf S}_{m\delta}(\tau)$.\\
Such algorithm is reported and proved in \cite{SPES2}; later Mora in a series of lecture notes labelled the restriction of this decomposition and interpolation to the case of a set of distinct points as `Axis-of-Evil' theorem and gave a precise description,  but no simple proof, of the result stated in \cite{SPES2}; S.  Steidel implemented the procedure in Singular \cite{DGPS}, \cite{St}.\\ 
We give here  such explicit algorithm that, given a finite set $\mathbf{X}$ of distinct points, provides a complete decomposition  $\mathbf{X} = \bigsqcup_{m=1}^n \bigsqcup_{\delta=1}^{d_m} {\sf S}_{m\delta}(\tau)$  on which, applying \cemu algorithm and interpolation, produces the required linear factorization for a lexicographical minimal \gb $F=\{f_1,...,f_r\}$ of the ideal $I(\mathbf{X})$ and thus a very simple proof of the `Axis-of-Evil' theorem in this particular situation.\\
This algorithm  arranges the $r$ terms $t_i$ belonging to $G(I(\mathbf{X}))$ with respect to lex ($t_1\leq...\leq t_r$) and constructs the factorization of each $f_i \in F$ through a suitable interpolation on a subset ${\sf S}_{m\delta}(t_i)$ of $\mathbf{X}$ depending on the exponents of the corresponding $t_i$.
More precisely, \cemu give an algorithm that enables to find the Groebner escalier $N(I(\mathbf{X}))$ and the minimal basis $G(I(\mathbf{X}))$ of the monomial ideal $T(I(\mathbf{X}))$.\\
Since the `Axis of Evil' algorithm's starting point are the elements of $\mathbf{X}$ and the  monomials of the finite Groebner escalier $N$ (computed using \cemu algorithm), but the algorithm requires as input the monomial basis of $T(I(\mathbf{X}))$, we also define the `potential expansion' 's algorithm.\\
It takes $N$ and computes the minimal basis.
\\
I note here that Marinari-Mora explicitly deduced,  as  trivial corollaries of  their `Axis-of-Evil' procedure, Lazard theorem (\cite{La}), Elimination theorem (\cite{B}),  Kalkbrener theorem (\cite{K2}), part of Gianni-Kalkbrener theorem (\cite{G},\cite{K}); they however remarked that, having being strongly influenced by  Gianni-Kalkbrner result, they cannot dismiss the possibility that  Gianni-Kalkbrenner argument is an essential tool of their proof of the `Axis-of-Evil' theorem.
\section{Notation.}\label{notation}
Let $P:=k[x_1,...,x_n]=\bigoplus_{d \in \NN}P_d$ be the ring of polynomials in $n$ variables and coefficients in the base field $k$. For all $M\subseteq P$, $M_d=M\cup P_d$ is its degree $d$ part.
Call $\mathcal{T}$ the semigroup of terms, generated by the set $\{x_1,...,x_n\}$:
$$\mathcal{T}:=\{x_1^{a_1}\cdots x_n^{a_n},\,(a_1,...,a_n)\in \NN^n \}.$$
Letting $\alpha=(\alpha_1,...,\alpha_n) \in \NN^n$, we will often write $x^{\alpha}$ instead of $x_1^{\alpha_1}\cdots x_n^{\alpha_n}$.\\
Define also the set 
$$T[m]:=\mathcal{T}\cap k[x_1,...,x_m]=\{x_1^{a_1}\cdots x_m^{a_m}/\, (a_1,...,a_m)\in \NN^m \}.$$ 
For each \emph{semigroup ordering} $<$ on $\mathcal{T}$ (i.e. a total ordering \st 
$ t_1<t_2 \Rightarrow tt_1<tt_2,\, \forall t,t_1,t_2 \in \mathcal{T}$) we can represent a polynomial $f\in P$ as a linear combination (with coefficients in $k$) of monomials arranged w.r.t. $<$:
$$f=\sum_{t \in \mathcal{T}}c(f,t)t=\sum_{i=1}^s c(f,t_i)t_i:\, c(f,t_i)\in k^*,\, t_i\in \mathcal{T},\, t_1>...>t_s.$$
We will call $T(f)=Lt(f):=t_1$  the \emph{leading term} of $f$ and $tail(f)=f-T(f)$ the \emph{tail} of $f$.
\\
We can also express it in a unique way as 
$$f=\sum_{i=0}^{\delta} g_ix_n^i\in k[x_1,...,x_{n-1}][x_n],\, g_i \in  k[x_1,...,x_{n-1}], \, g_{\delta}\neq 0$$
(where $\delta:=deg_n(f)$ is the degree w.r.t. $x_n$).\\
We denote $Lp(f):=g_{\delta}$, the \emph{leading polynomial} of $f$.
\begin{definition}\label{Predecessori}
For each monomial $t \in \mathcal{T}$ and $x_j \vert  t$, the only $u \in \mathcal{T}$ \st $t=x_ju$ is called $j$-th \emph{predecessor} of $t$. 
\end{definition}
A subset $N\subseteq \mathcal{T}$ is an \emph{order ideal} if 
$$t \in N \Rightarrow s \in N\, \forall s \vert t .$$
Let $N \subset \mathcal{T}$ an order ideal, 
A subset $N\subseteq \mathcal{T}$ is an order ideal if and only if $\mathcal{T}\setminus N=J$ is a semigroup ideal (i.e. $\tau \in J \Rightarrow t\tau \in J,\, \forall t \in \mathcal{T}$).\\
We set $N(J):=\mathcal{T}\setminus T(J)=N$.
\\
For a semigroup ideal $J$, $G(J)$ denotes its minimal basis and
\begin{eqnarray*}
&G(J):=\{\tau \in J\, \vert \, \textrm{ each predecessor of }\, \tau \in N(J)\}=&\\
&=\{\tau \in \mathcal{T}\, \vert \,N(J)\cup\{\tau\}\, \textrm{order ideal}, \, \tau \notin N(J)\}.
\end{eqnarray*}
For all subsets $G \subset P$, we define $T\{G\}:=\{T(g),\, g \in G\}$ and we call $T(G)$ the semigroup ideal $\{\tau T(g),\, \tau \in \mathcal{T}, g \in G\}$, generated by $T\{G\}$.\\
For any ideal  $I \triangleleft P$ consider the semigroup ideal $T(I)=T\{I\}$, denoting by abuse of notation $G(I)$ its minimal basis $G(I)$ and the border ideal of $I$
\begin{eqnarray*}
&B(I):=\{x_ht,\, 1 \leq h \leq n,\, t \in N(I)\}\setminus N(I)=&\\
&=T(I)\cap (\{1\}\cup \{x_ht,\, 1 \leq h \leq n,\, t \in N(I)\}). 
\end{eqnarray*}
We will always consider the \emph{lexicographic order} induced by $x_1<...<x_n$, i.e:
$$x_1^{a_1}\cdots x_n^{a_n}<x_1^{b_1}\cdots x_n^{b_n} \Leftrightarrow \exists j\, \vert  \, a_j<b_j,\,a_i=b_i,\, \forall i>j. $$
This is a  \emph{term order}, that is a semigroup ordering \st $1$ lower to every variable or, equivalently, it is a \emph{well ordering}.
\begin{lemmadef}\label{lem1.1 IRAN} We have:
\begin{enumerate}
\item $P \cong I \oplus k[N(I)];$ 
\item $P/I \cong k[N(I)];$
\item $\forall f \in P$, $\exists ! g:=Can(f,I)=\sum_{t \in N(I)} \gamma(f,t,<)t \in k[N(I)]$, called \emph{canonical form} of $f$ with respect to $I$, \st $f-g \in I.$
\end{enumerate}
\end{lemmadef}
\begin{definition}\label{def1.2 IRAN}
Given a term order $ <$ on $\mathcal{T}$: 
\begin{enumerate}
	\item a \emph{\gb} of $I$ is a set $G \subset I$ \st  $T(G)=T\{I\}$, that is $T\{G\}$ generates the semigroup ideal  $T(I)=T\{I\}$;
	\item a \emph{minimal \gb } is a \gb such that divisibility relations among the leading monomials of its members do not exist;
	\item the unique \emph{reduced \gb} of $I$  is the set:
	$$\mathcal{G}(I):=\{\tau - Can(\tau, I):\, \tau \in G(I)\}.$$
	Each member of the reduced \gb has a leading term which does not divide any monomial of another member.
\end{enumerate}
\end{definition}
Let $\mathbf{X}=\{P_1,...,P_N\} \subset k^n$ be a finite set of distinct points 
$$P_i:=(a_{i1},...,a_{in}),\, i=1,...,N.$$
We call
$$I(\mathbf{X}):=\{f \in P:\, f(P_i)=0,\, \forall i\},$$
the \emph{ideal of points} of $\mathbf{X}$.
\\
Finally, we define  the projection maps:\\
\begin{minipage}[b]{0.5\linewidth}
\centering
$$\pi_m:k^n \rightarrow k^m $$
$$(X_1,..,X_n)\mapsto (X_1,...,X_m), $$
\end{minipage}
\hspace{0.05cm}
\begin{minipage}[b]{0.5\linewidth}
\centering
$$\pi^m:k^n \rightarrow k^{n-m+1} $$
$$(X_1,..,X_n)\mapsto (X_m,...,X_n)$$
\end{minipage}
and, for $P \in k^n,\; \mathbf{X} \subset k^n$, let
$$\Pi_s(P,\mathbf{X}):=\{P_i \in \mathbf{X}  / \pi_s(P_i)=\pi_s(P)\}, $$
$$\Pi^s(P,\mathbf{X}):=\{P_i \in \mathbf{X}  / \pi^s(P_i)=\pi^s(P)\}, $$
extending in the obvious way the meanings of $\pi_s(\mathbf{d}),\pi^s(\mathbf{d}),\Pi_s(\mathbf{d},D), \Pi^s(\mathbf{d},D)$ to $\mathbf{d}\in \NN^n\subset k^n$ e $D \subset \NN^n \subseteq \NN^n.$
\\
With the same notation $\pi_m$ we  indicate also
$$\pi_m: \mathcal{T}\cong \NN^n \rightarrow \NN^m\cong T[m] $$
$$x_1^{a_1}\cdots x_n^{a_n}\mapsto x_1^{a_1}\cdots x_m^{a_m}.$$
  
\section{\cemu Correspondence.}\label{CEMU}
Cerlienco and Mureddu (\cite{CeMu}, \cite{CeMu2}, \cite{CeMu3}) provided an algorithm which solves the following\\
\begin{prob}\label{5.1 IRAN}
Given  finite ordered set of distinct points 
$$\underline{\mathbf{X}}:=(P_1,...,P_N)\subset k^n;\; P_i:=(a_{i1},...,a_{in}) $$
compute a monomial basis (w.r.t. the lexicographic order induced by $x_1<...<x_n$) of the quotient $k[x_1,...,x_n]/I(\mathbf{X})$, where  $\mathbf{X}$ denotes the support $\{P_1,...,P_N\}$ of $\mathbf{X}$.
\end{prob}\\
More precisely, they 
\begin{itemize}
\item define the operator $\Phi$, associating to  $\underline{\mathbf{X}}$ an ordered set $$\Phi(\underline{\mathbf{X}}):=(\mathbf{d}_1,...,\mathbf{d}_N)\subset \NN^n$$ \st 
$\vert\Phi(\underline{\mathbf{X}})\vert=\vert \underline{\mathbf{X}}\vert=N $
and \st, for all $m<N$ the subset $(\mathbf{d}_1,..\mathbf{d}_m)$ is exactly
$\Phi((P_1,...,P_m) )$.
\item define the $\mathbf{\sigma}$\emph{-value}  w.r.t. $\mathbf{X}$ $s=\sigma(P,\mathbf{X})$  of a point $P \in K^n\setminus \mathbf{X}$ as the maximal integer \st $\Pi_{s-1}(P,\mathbf{X})\neq \emptyset$ (by convention, $\forall P,\mathbf{X}$, $\Pi_0(P,\mathbf{X})\neq \emptyset$).\\
For $P \notin \mathbf{X}$, they define the set 
$$\Sigma(P,\mathbf{X}):=\{P_i \in \mathbf{X}/ \pi_{s-1}(P_i)=\pi_{s-1}(P),\, s=\sigma(P,\mathbf{X})\} $$
containing all the points of $\mathbf{X}$ having the first $s-1$ coordinates equal to those  of $P\notin \mathbf{X}$. They extend the notation to the case $P=P_j \in \underline{\mathbf{X}}$ in the following way:
$$\sigma(P, \underline{\mathbf{X}}):=\sigma(P,\{P_1,..,P_{j-1}\})$$
$$\Sigma(P, \underline{\mathbf{X}}):=\Sigma(P,\{P_1,..,P_{j-1}\}).$$
\end{itemize}
\begin{oss}
Given a term order $\preceq$, a monomial basis for $A:=k[x_1,...,x_n]/I(\mathbf{X})$,$$[\mathbf{x}^{\mathbf{i_1}}],...,[\mathbf{x}^{\mathbf{i_N}}], \; with \; \mathbf{x}^{\mathbf{i_1}}\preceq...\preceq \mathbf{x}^{\mathbf{i_N}}$$
 is called \emph{minimal} with respect to the term order if, 
for every other monomial basis $[\mathbf{x}^{\mathbf{i_1'}}],...,[\mathbf{x}^{\mathbf{i_N'}}]$, with $\mathbf{x}^{\mathbf{i_1'}}\preceq...\preceq \mathbf{x}^{\mathbf{i_N'}}$ for the $A$
 it holds 
$$\forall j=1,...,N,\; \mathbf{x}^{\mathbf{i_j}}\preceq \mathbf{x}^{\mathbf{i_j'}}.$$ In \cite{CeMu}, they state that the computed monomial basis is the minimal one.
 \end{oss}
\begin{algorithm}
\caption{\cemu algorithm.}
\begin{algorithmic}[1]
\Procedure{CeMu}{$\underline{\mathbf{X}}$} $\rightarrow \Phi(\underline{\mathbf{X}})$
\If{$N=1$}
\State $\Phi(\underline{\mathbf{X}}):=\{(0,...,0)\}$.
\EndIf
\If{$1 < N$}\Comment{{\tiny suppose to know by induction hypothesis $\Phi((P_1,...,P_{N-1}))=(\mathbf{d}_1,...,\mathbf{d}_{N-1})$ and look for $\mathbf{d}_N=\Phi(P_N)$.}}
\State $s=\sigma(P_N,\underline{\mathbf{X}})$.
\For{$i=n$ to $1$}
\If{$i>s$}
\State  $d_{Ni}=0$.
\EndIf
\If{$i=s$}
\State $m,\,(1 \leq m \leq n)$, maximal s.t $\,\pi_{s-1}(P_m)=\pi_{s-1}(P_N)$, $\pi^{s+1}(\mathbf{d}_m)=(0,...,0)=\pi^{s+1}(\mathbf{d}_N)$.\Comment{{\tiny  $P_m$ is the  $\mathbf{\sigma}$~\emph{-antecedent} of $P_N$ w.r.t. $(P_1,...,P_{N-1}),$  $\Phi((P_1,...,P_{N-1})).$}}
\State $d_{Ns}=d_{ms}+1$.
\EndIf
\If{$i<s$}\Comment{{\tiny we use induction here.}}
\State $\mathcal{W}(P_N,\underline{\mathbf{X}}):=\{P \in \underline{\mathbf{X}} \vert \,\Phi(P)=\mathbf{d}=(*,...*,d_{Ns},0,...,0),\,\}=\{P_{j1},...,P_{jr}\}.$
\State $\mathcal{Q}:=\pi_{s-1}(\mathcal{W}(P_N,\underline{\mathbf{X}})).$\Comment{{\tiny $\vert\mathcal{Q}\vert=\vert\mathcal{W}(P_N,\underline{\mathbf{X}})\vert=r<N$. If $h<r=\vert\mathcal{W}(P_N,\underline{\mathbf{X}})\vert$, then $\pi_{s-1}(P_{jh})\neq \pi_{s-1}(P_N).$
Moreover, since $\Phi$ is inductive, if $h<k\leq r$ then $\pi_{s-1}(P_{jh})\neq \pi_{s-1}(P_{jk})$.}}
\State $\pi_{s-1}(\mathbf{d}_N)= \widetilde{\mathbf{d}_r}.$\Comment{{\tiny By the induction hypothesis,$\Phi(\underline{\mathcal{Q}})=(\widetilde{\mathbf{d}_1},.., \widetilde{\mathbf{d}_r})$ and $\forall 1 \leq i < r$, $\widetilde{\mathbf{d}_i}=\pi_{s-1}(\mathbf{d}_{ji})$.}}
\State \textbf{break.}
\EndIf
\EndFor
\EndIf\\
\Return $\Phi(\underline{\mathbf{X}}).$
\EndProcedure
\end{algorithmic}
\end{algorithm}
\begin{proposition}({\cite{CeMu}})\label{prop6 CEMU}\\
Let $D:=\Phi(\mathbf{X})$. Then $\{[\mathbf{x^d}]/ \mathbf{d}\in D\}$ is a monomial basis for $k[x_1,...,x_n]/I(\mathbf{X}).$
\\
Such a monomial basis is minimal with respect to the given $<$.
\end{proposition}
Once the Groebner escalier $N$ is known, it is very simple to compute the minimal basis $G$ of $T(I(\mathbf{X}))=\mathcal{T}\setminus N$.
Given the set $\mathbf{X}$, the first step to compute the linear factorization of a minimal \gb will be to apply \cemu algorithm to $\mathbf{X}$ and compute $N$, in order to obtain $G$.
\section{The potential expansion's algorithm.}
Consider the polynomial ring $k[x_1,...,x_n]$ with usual ordering  $<$. Given a finite set of distinct points $\mathbf{X}=\{P_1,...,P_N\}$, consider the ideal $I(\mathbf{X})\triangleleft k[x_1,...,x_n]$ which is radical and zerodimensional, so its Groebner escalier $N$ is a finite set.\\
We will compute the minimal monomial basis  $G$ of $T(I(\mathbf{X}))$, given the Groebner escalier. The algorithm actually provides correct results irrispective of the given term ordering, but 
since we use \cemu correspondence, we will consider only our lex order.\\
In order to make the reasoning clear, we will represent the monomials using the same diagrams introduced in \cite{MR} to study properties of Borel ideals.
\\
Apply \cemu correspondence to $\mathbf{X}$ in order to have $N(\mathbf{X})=\{\tau_1,...,\tau_N\}$. It is well known (see, for instance \cite{SPES2}) that $\vert N(\mathbf{X}) \vert =\vert \mathbf{X} \vert$.\\
We first define the potential expansion of a subset $H \subset \mathcal{T}$, from which the algortihm  bears its name.
\begin{definition}\label{espansione potenziale}
Let $H \subseteq\mathcal{T}_j$ for some $j \in \NN^*$ we set $C^{(0)}(H):=H$ and, for all $l \in \NN^*$ $C^{(l)}(\tau)=\mathcal{T}_{j+l}\setminus\{x_1,...,x_n\}\cdot (\mathcal{T}_{j+l-1}\setminus C^{(l-1)}(H))$.\end{definition}
We then slice the Groebner escalier by degree, having $N_0,\,N_1,\, \cdots \, N_h,$ where $h$ is the maximal degree of terms appearing in $N$.\\
The minimal monomial basis $G(I(\mathbf{X}))$ will have at most degree $h+1$.
As a matter of fact, if $\tau \in G$ with $deg(\tau)=d>h+1$ its predecessors will belong to $N$ and have degree $d-1\geq h+1$ which is impossible.\\
The computation of $G$ is performed as follows. \\
Consider $\mathcal{T}_i$ $\forall i=0,...,h+1$: it is well known that $\vert \mathcal{T}_i\vert ={n+i-1\choose n-1}$.\\
For each $i$, define $Gen_i(I):=\{t \in G(I)\vert \, deg(t)\leq i\}$. Since $I$ is a proper ideal, $Gen_0(I)=\emptyset$.
\\
Let $h$ the minimal $i$ such that $Gen_h(I) \neq \emptyset$, $\forall i\geq 1$
$$Gen_{i+h}=Gen_{h+i-1}\cup(\mathcal{T}_{h+i}\setminus (N_{h+i}\cup \bigcup_{j=h+1}^{h+i-1} C^{(h+i-j)}(G_j))).$$
We then have
\begin{algorithm}
\caption{The potential expansion's algorithm.}
\begin{algorithmic}[1]
\Procedure{PotExp}{$N(I)$} $\rightarrow I$\Comment{{\tiny $I$ is expressed using its minimal basis.}}
\Require  $N=[N_0,...,N_h,N_{h+1}]$, \st $N_{h+1}=\emptyset$.
\State $C=[\emptyset]$.\Comment{{\tiny the potential expansion' s list.}}
\State $Gen=\emptyset$.
\State $I=(0)$.
\For{$i=0$ to $h+1$}
\State  $d={n+deg( N_i[1])-1\choose n-1}-\vert N_i\cup C[i] \vert$.
\If{$d=0$} \Comment{{\tiny no new generators.}}
\State  $C[i+1]=PotentialExpansion(C[i])$.
\State $Gen_i=(0)$
\Else \Comment{{\tiny adding new generators.}}
\State $Gen_i=\mathcal{T}_i\setminus (N_i\cup C[i])$.
\State  $C[i+1]=PotentialExpansion(Gen_i\cup C[i])$.
\State  $I=I+Gen_i$.
\EndIf
\EndFor\\
\Return $I$
\EndProcedure
\end{algorithmic}
\end{algorithm}\\
The algorithm uses a subroutine $PotentialExpansion$ \st
$$PotentialExpansion(F)=C^{(1)}(F).$$
We will also have a subroutine finding $\mathcal{T}_{h+i}\setminus (N_{h+i}\cup \bigcup_{j=h+1}^{h+i-1} C^{(h+i-j)}(G_j))$.\\
WLOG we will think that the sets $\mathcal{T}_{h+i}$ and $N_{h+i}\cup \bigcup_{j=h+1}^{h+i-1} C^{(h+i-j)}(G_j)$ are ordered w.r.t. the same ordering, since it is enough  to perform a merging with the Groebner escalier and the \pe previously ordered.\\
All these steps end: the subroutine finding the complementary can be developed performing a loop on the two ordered lists $A:=\mathcal{T}_i=[a_1,..,a_m], \, m\geq n$ and $B:=N_i\cup C^{(i)}=[b_1,...,b_n]$ (using two indices $i,j$), keeping in mind that $B \subseteq A$ or $B=A$ and that  $B[j]\geq A[i]$ at every step. Start with $b_1$: if $b_1=a_1$ we set $i=i+1; j=j+1$.\\
If we find $a_i\neq b_j$ for a certain couple $(i,j)$, we put  $A[i]$ in the complementary and $i=i+1$ without modifying $j$.
\begin{example}
There are situations in which  $N$ contains monomials of degree at most $h$, but also the minimal basis shares the same property.
\\
Take $I=(x^3,y^2,z^2,xy)\triangleleft k[x,y,z]$, whose Groebner escalier is:\\
$N_0=\{1\}$\\	
$N_1=\{x,y,z\}$\\
$N_2=\{yz,xz,x^2\}$\\
$N_3=\{x^2z\}$:
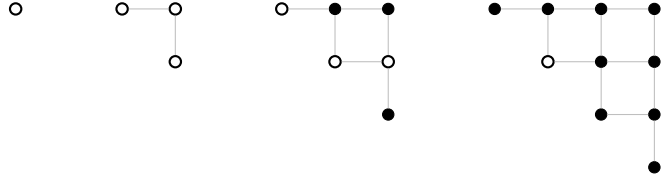
\begin{center}
\begin{tikzpicture}[scale=0.35]
\tikzstyle{ideal}=[circle,draw=black,fill=black,inner sep=1.5pt]
\tikzstyle{quotient}=[circle,draw=black,thick,inner sep=1.5pt]
\tikzstyle{basis}=[circle,draw=red,fill=red,inner sep=1.5pt]
\node (1) at (-12,0) [quotient] {};
\node (2) at (-8,0) [quotient] {};
\node (3) at (-6,0) [quotient] {};
\node (4) at (-2,0) [quotient] {};
\node (5) at (0,0) [ideal] {};
\node (6) at (2,0) [ideal] {};
\node (7) at (-6,-2) [quotient] {};
\node (8) at (0,-2) [quotient] {};
\node (9) at (2,-2) [quotient] {};
\node (10) at (2,-4) [ideal] {};
\node (11) at (6,0) [ideal] {};
\node (12) at (8,0) [ideal] {};
\node (13) at (10,0) [ideal] {};
\node (14) at (12,0) [ideal] {};
\node (15) at (8,-2) [quotient] {};
\node (16) at (10,-2) [ideal] {};
\node (17) at (12,-2) [ideal] {};
\node (18) at (10,-4) [ideal] {};
\node (19) at (12,-4) [ideal] {};
\node (20) at (12,-6) [ideal] {};
\draw [-,black!25] (2) -- (3);
\draw [-,black!25] (4) -- (5);
\draw [-,black!25] (3) -- (7);
\draw [-,black!25] (5) -- (8);
\draw [-,black!25] (5) -- (6);
\draw [-,black!25] (6) -- (9);
\draw [-,black!25] (9) -- (10);
\draw [-,black!25] (8) -- (9);
\draw [-,black!25] (11) -- (12);
\draw [-,black!25] (15) -- (12);
\draw [-,black!25] (13) -- (12);
\draw [-,black!25] (13) -- (14);
\draw [-,black!25] (13) -- (16);
\draw [-,black!25] (14) -- (17);
\draw [-,black!25] (15) -- (16);
\draw [-,black!25] (16) -- (17);
\draw [-,black!25] (16) -- (18);
\draw [-,black!25] (17) -- (19);
\draw [-,black!25] (18) -- (19);
\draw [-,black!25] (20) -- (19);
\end{tikzpicture}
\end{center}
The monomial basis does not contain elements of degree $4$.
\end{example}
We call $G_i$ the set of $i$-degree elements  of the minimal basis and $I$ the monomial ideal we want to find.\\
\begin{lemma}\label{NuC complem G}
For all $i=0,...,h+1$
$$T_i\setminus (N_i\cup \bigcup_{j=1}^{i-1}C^{i-j}(G_j) )=G_i.$$
\end{lemma}
\begin{dimo}
The inclusion $T_i\setminus (N_i\cup \bigcup_{j=1}^{i-1}C^{i-j}(G_j) )\supseteq G_i$ is trivial, so we only prove  $T_i\setminus (N_i\cup \bigcup_{j=1}^{i-1}C^{i-j}(G_j) )\subseteq G_i$.\\
Consider $\tau \in \mathcal{T}_i\setminus  (N_i\cup \bigcup_{j=1}^{i-1}C^{i-j}(G_j) )$. Clearly $\tau \in I$. 
\\
Let $\sigma$ the $i$th predecessor of $\tau$; if  $\sigma \in I$, $\exists t \in G$ with $\sigma =t\cdot m$ for a suitable $m \in \mathcal{T}$.\\
Then $\tau =t\cdot m \cdot x_i$ i.e. $\tau \in \bigcup_{j=1}^{i-1}C^{i-j}(G_j)$.
\end{dimo}\\
This lemma assures that the result obtained via the potential expansion's algorithm is correct.
\section{The Axis of Evil Algorithm.}\label{AoE}
A $\mathbf{0-}$\emph{dimensional radical} ideal $I \triangleleft P$ is completely determined if we know the set $V(I)$ of its zeros.
\\
Consider a finite set of distinct points $\mathbf{X}=\{P_1,...,P_r\}$; we will denote indifferently the Groebner escalier of the ideal $I(\mathbf{X})$ with $N(I(\mathbf{X}))$ or $N$. A variation of Cerlienco-Mureddu algorithm (\cite{CeMu}) allows us to find a `linear factorization' for every element of a lexicographic minimal \gb in the sense of the
\begin{theorem}\label{ASSE DEL MALE}
Let $t_i:=x_1^{d_1}\cdots x_n^{d_n},\, i=1,...,r$ be the generators of the minimal basis of $T(I)$, where $I$ is a $0-$ dimensional radical ideal.
\\
A combinatorical algorithm and interpolation allow us to deduce polynomials
$$\gamma_{m\delta i}=x_m-g_{m\delta i}(x_1,...,x_{m-1}), $$
$\forall i,m,\delta$, with $1 \leq i \leq r$,  $1 \leq m \leq n$,  $1 \leq \delta \leq d_m$ \st
$$f_i= \prod_{m} \prod_{\delta}\gamma_{m\delta i}\; \forall i$$
where $f_i,\, i=1,...,r$ are the polynomials forming a minimal \gb of $I$ with respect to the lexicographic order induced by $x_1<...<x_n$.
\end{theorem}
In that algorithm we will use the projections, as we defined in section \ref{CEMU}.
\begin{algorithm}\label{ASSEMALE}
\caption{The Axis of Evil algorithm.}

\begin{algorithmic}[1]
\Procedure{AoE}{$\mathbf{X},G(I(\mathbf{X})):=\{\tau_1,...,\tau_r\}$} $\rightarrow R$\Comment{{\tiny $R$ contains a factorized minimal \gb of $I$.}}
\Require the elements $G(I(\mathbf{X}))$ are in increasing order w.r.t the lexicographical order w.r.t. $x_1<...<x_r$.
\State $R=\emptyset$
\For{$i=1$ to $r$}
\State $N_1(\tau_j):=\{x_1^i/\; i< d_1\}=\{\omega \in T[1],\, \tau_j > \omega x_{2}^{d_{2}} \cdots x_n^{d_n}\in N\}$
\State $A_1(\tau_j):=\{\Phi^{-1}(x_1^i x_2^{d_2}\cdots x_n^{d_n})/\, i<d_i \}\subset \mathbf{X}$.\
\State $B_1(\tau_j):=\pi_1(A_1(\tau_j))\subset k$.
\State $\gamma_{1\tau_j}:=\prod_{a \in B_1(\tau_j)}(x_1-a)$.
\For{$m=2$ to $n$}
\State  $\zeta_{m \tau_j} := \prod_{\nu =1}^{m-1}\gamma_{\nu \tau_j}$.
\State $D_{m0}:=\{P_i \in \mathbf{X}/\, \zeta_{m\tau_j}(P_i) \neq 0\}$.
\If{$\vert D_{m0}\vert=0$}
\State $R=[R,\zeta_{m\tau_j}]$.
\State \textbf{break.}
\EndIf
\State $N_m(\tau_j):=\{\omega \in T[m],\, \tau_j > \omega x_{m+1}^{d_{m+1}} \cdots x_n^{d_n}\in N \}$.
\For{$\delta=1$ to $d_m$}
\State $A_{m\delta}(\tau_j):=\{\Phi^{-1}(vx_m^{d_m-\delta}x_{m+1}^{d_{m+1}}\cdots x_n^{d_n})/\, v \in T[m-1],\, vx_m^{d_m-\delta} \in N_m(\tau_j)\} \cap D_{m(\delta -1)}(\tau_j)$.
\State $E_{m\delta}(\tau_j):=\Phi(\pi_m(A_{m\delta}(\tau_j)))$.
\State      \begin{equation*}
				                               \gamma_{m\delta \tau_j}:=x_m+\sum_{\omega \in E_{m\delta}(\tau_j)}                    				                               c(\gamma_{m\tau_j}, \omega)\omega,
				                                \label{eq:interpolyn }
				                                \end{equation*}
				                                \st $\gamma_{m\delta \tau_j}(P)=0,\, \forall P \in A_{m\delta}(\tau_j).$ 
	\State $\xi_{m \delta}:=\prod_{\nu=1}^{m-1}\gamma_{\nu \tau_j}\prod_{d=1}^{\delta} \gamma_{m d \tau}.$           
	\State $D_{m\delta }(\tau_j):=\{P_i \in \mathbf{X}/\, \xi_{m\delta}(P_i)\neq 0\} \subseteq \mathbf{X}$
	\If{$\vert D_{m\delta}(\tau_j)\vert=0$}
	\State $R=[R,\xi_{m\delta}]$.
	\State \textbf{break.}
	\EndIf
\EndFor
\State $\gamma_{m\tau_j}:=\prod_{\delta}\gamma_{m\delta \tau_j}$.
\EndFor
\EndFor\\
\Return $R.$
\EndProcedure
\end{algorithmic}
\end{algorithm}\\
The Axis of Evil algorithm works then in the following way:
\begin{itemize}
\item consider $\tau_j :=x_1^{d_1}\cdots x_n^{d_n} \in G$. The required polynomial $f=\tau_j+tail(f)$ is factorized in $\sum_{i=1}^n d_i$ factors: $d_1$ polynomials whose leading term is $x_1$, $d_2$ polynomials \st their leading term is $x_2$ and so on;
\item consider the monomials $x_1^{a_1}x_2^{d_2}\cdots x_n^{d_n}$ such that $a_1<d_1$;
\item every such monomial is associated, via \cemu Correspondence, to a point of our set $\mathbf{X}$. Project these points with respect to the first coordinate, obtaining $d_1$ numbers $y_1,...,y_{d_1}$;
\item $x_1-y_i$, $i=1,...,d_1$ are the first $d_1$ factors;
\item construct the subset $D_{20}$ of $\mathbf{X}$ containing all the points in which the product $(x_1-y_1)\cdots(x_1-y_{d_1})$ does not vanish. If it is empty then stop and consider the next monomial in $G$; otherwise continue as follows;
\item find the set $N_2(\tau_j)$ of all monomials in $\mathcal{T}[2]$ \st $x_1^{\alpha_1}x_2^{\alpha_2}<x_1^{d_1}x_2^{d_2}$; 
\item split the elements of $N_2(\tau_j)$ with respect to the exponents of $x_2$ and construct, via \cemu correspondence, the set 
$$\{\Phi^{-1}(vx_2^{d_2-\delta}x_{3}^{d_{3}}\cdots x_n^{d_n})/\, v \in T[1],\, vx_2^{d_2-\delta} \in N_2(\tau_j)\}$$
\item intersect the previous set with $D_{20}$, project the resulting set of points ($A_{2\delta}(\tau_j)$) with respect to the first two coordinates and apply \cemu Correspondence, obtaining a set $E_{2\delta \tau_j}$;
\item interpolate over $A_{2\delta}(\tau_j)$, finding $d_2$ factors whose leading terms are all equal to $x_2$. The monomials of $E_{2\delta \tau_j}$  are the ones appearing in such factorization;
\item update the set of points in which the current polynomial does not vanish and stop if it is empty;
\item repeat these steps letting all the variables vary one by one; 
\item repeat all the steps for all $\tau_i \in G$.
\end{itemize}
\begin{oss}
{\rm Given $\tau_j=x_1^{d_1}\cdots x_n^{d_n} \in G$, every variable $x_i$ will appear only $d_i$ times in the execution of the algorithm.}
\end{oss}
\begin{oss}
{\rm The sets $N_m(\tau_j):=\{\omega \in T[m],\, \tau_j > \omega x_{m+1}^{d_{m+1}} \cdots x_n^{d_n}\in N \}$ (in particular for $m=1$ we have $N_1(\tau_j):=\{x_1^i/\; i< d_1\}$)  are constructed in order to determine in which points  it is necessary to interpolate.\\
Since for $\mu >\tau_j$ the \cemu correspondence provides a point $P_{\mu '}$ such that $\exists k \in \{1,...,n\}:\; \pi_k(P_{\mu})=\pi_k(P_{\mu '})$, in order to obtain polynomials vanishing on all the poinst of $\mathbf{X}$ it is not necessary to interpolate in the whole $\Phi^{-1}(N)$ as it suffices to consider only those corresponding to $\mu \in N$ with $\mu <\tau_j$.}
\end{oss}
\begin{oss}
{\rm The terms smaller than $\tau_j$ mentioned before are found 
releasing all the variables one by one.\\
Imagine the monomials in $k[x_1,...x_n]$ as points in $k^n$, identifying every term to the $n$-uple of its exponents. So we can `draw' them in a $n$-dimensional space and we can think our realeasings as an increment by one of the `directions' where we can move there.\\
We point out that $N_m(\tau_j)\subseteq N_h(\tau_j)$ for $m\leq h$. \\If $\omega \in N_m(\tau_j)$, 
$\tau_j>\omega x_{m+1}^{d_{m+1}}\cdots x_n^{d_n}\in N$; as $\omega x_{h+1}^{d_{h+1}}\cdots x_n^{d_n} \vert x_{m+1}^{d_{m+1}}\cdots x_n^{d_n}$ we have $\omega x_{h+1}^{d_{h+1}}\cdots x_n^{d_n}\in N$ and\\ 
$\omega x_{h+1}^{d_{h+1}}\cdots x_n^{d_n} \leq x_{m+1}^{d_{m+1}}\cdots x_n^{d_n}<\tau_j$.\\
At each step we will count out all the points in which the polynomial already vanishes and we will stop the computation when the current factorized polynomial vanishes on the whole $\mathbf{X}$.\\
We will see an example of it later.}
\end{oss}
\begin{oss}\label{delta0}
{\rm If the number of released variables is $>1$, we also must split the obtained monomials regarding the exponent of the maximal variable.\\
Consider then the loop on $\delta$ and, in particular, the set:
 $$C_{m\delta}(\tau_j):=\{\Phi^{-1}(vx_m^{d_m-\delta}x_{m+1}^{d_{m+1}}\cdots x_n^{d_n})/\, v \in T[m-1],\, vx_m^{d_m-\delta} \in N_m(\tau_j)\}.$$
We intersect $C_{m\delta}(\tau_j)$  with the subset of $\mathbf{X}$ containing the points not vanishing the  current factorized polynomial.\\
 We can easily notice that, performing the algorithm, we only compute the sets $C_{m1}(\tau_j),...,C_{md_m}(\tau_j)$, but in  $N_m(\tau_j)$ there are also monomials   $\omega=x^{a_1}_1\cdots x^{a_{m-1}}_{m-1}x_m^{d_m}$ \st $\tau_j > \omega x_{m+1}^{d_{m+1}} \cdots x_n^{d_n}\in N $, which would be generated considering $\delta =0$.\\
They   are not considered in the algorithm  because they are related to monomials examined in the previous step:
$=x_1^{a_1}\cdots x_{m-1}^{a_{m-1}}\in N_{m-1}$, so the corresponding points have already been treated. Taking  $\delta=0,..,d_m$, the sets $C_{m\delta}(\tau_j)$ form a partition of $N_m(\tau_j)$ basing on the degree of $x_m$. As a matter of fact, in order to have  $\omega \in N_m(\tau_j)$ we must have $\tau_j > \omega x_{m+1}^{d_{m+1}} \cdots x_n^{d_n}$, where $\omega x_{m+1}^{d_{m+1}} \cdots x_n^{d_n} \in N,$  then the exponent of $x_m$ will be the first checked in the lexicographic test and so it will be limited by $d_m.$\\ 
According to the values of this exponent, the ones associated to smaller variables will vary.}
\end{oss}
\begin{oss}\label{OrdinePt}
{\rm At the beginning of the algorithm, we imposed the monomials  $\tau_j,$\\
$ j=1...,r$ to be in increasing order with respect  $<$.
The steps made by the algorithm on each $\tau_j$ are totally independent both on those made and on those to be made on a monomial $\tau_k$ (it is indifferent whether $j\gtrless k$) belonging to $G$, so we will obtain the same factorizations even if we launch the computation on a list of unordered monomials.
\\
Clearly, the result of our computation is not the reduced \gb of the given ideal, it is only one of the minimal \gbx but we can obtain the reduced \gb via simple reduction.
\\
We decided to put the monomials in such an order because we want \emph{every polynomial to be reduced with respect to the `previous' ones}.\\
If $f_j$  is one of our resulting polynomials and $Lt(f_j)=\tau_j$, the polynomials utilizable to reduce $f_j$ (the previous ones)  must be necessarily all and only the ones having as leading terms elements in $G$ lower than the given $\tau_j$.
}
\end{oss}
The algorithm terminates because it works on:
\begin{enumerate}
\item points in the finite set $\mathbf{X}$;
\item monomials $\tau \in G$ (they are in a finite number, \cite{SPES2});
\item a finite set of variables.
\end{enumerate}
Let us study the correctness of the algorithm.
\begin{lemma}\label{si annulla ovunque}
The factorized polynomials obtained from our algorithm  vanish on all the points of the set $\mathbf{X}$.
\end{lemma}
\begin{dimo}
Suppose we want to  construct $\gamma_{\tau}$ with $\tau =x_1^{\alpha_1}\cdots x_n^{\alpha_n}$. \\
Let $\mu =x_1^{\beta_1}\cdots x_n^{\beta_n}$, corresponding to a point $P_{\mu}\in \mathbf{X}$ through \cemu Correspondence.\\
Let $\mu < \tau$, then at least one of the exponents of the variables appearing in  $\mu$ is lower than the corresponding in $\tau$, say $\beta_i<\alpha_i$, so $\mu $ is linked to an element of  $ N_i(\tau)$ and so it can, alternatively: 
\begin{itemize}
\item belong to $A_{i\delta}(\tau)$ for some $\delta$;
\item be such that the corresponding point already annihilates the polynomial found.
\end{itemize}
If   $\mu > \tau$ (since $\tau \notin N$, it is surely impossible that $\tau =\mu$) then there will be a point $P_{\mu'}$ \st
$$\pi_j(P_{\mu})=\pi_j(P_{\mu'}), $$
corresponding to a  $\mu' < \tau$.\\
We then use $\mu'$ and we come back to the previous situation.
\end{dimo}
\begin{corollario}\label{corAoE}
The ideal generated by these polynomials is exactly $I(\mathbf{X})$. 
\end{corollario}
\begin{dimo}
By the previous lemma, the polynomials vanish on all the points of the set  $\mathbf{X}$ and the equality comes out by reasons of multiplicity
\end{dimo}\\
The resulting polynomials form a minimal \gb because:
\begin{itemize}
\item they vanish on all the points of $\mathbf{X}$;
\item their heads form exactly $G(I(\mathbf{X}))$.
\end{itemize}
 Notice that we can obtain the current interpolating polynomial applying Moeller algorithm to the projection through $\pi_m$ of the points of the current $A_{m\delta}(\tau)$ (\cite{MMM}). 
\begin{example}\label{Esempio Aoe2}
Let\\
$\mathbf{X}:=\{(4,0,0),(2,1,4),(2,4,0),(3,0,1),(2,1,3),(1,3,4),(2,4,3),(2,4,2),(1,0,2)\}.$\\
$P_1:=( 4,0,0):$ it is a single point, so $\Phi(\{(4,0,0)\})=(0,0,0)$\\
$P_2:=(2,1,4):$ $s=1$, $m=1$, $(1,0,0)$\\
$P_3:=(2,4,0):$ $s=2$, $m=2$,  $(0,1,0)$\\
$P_4:=(3,0,1):$ $s=1$, $m=1$,  $(2,0,0)$\\
$P_5:=(2,1,3):$ $s=3$, $m=2$,  $(0,0,1)$\\
$P_6:=(1,3,4):$ $s=1$, $m=4$,  $(3,0,0)$\\
$P_7:=(2,4,3):$ $s=3$, $m=3$, $\mathbf{ W}=\{(2,1,3),(2,4,3)\}$,   $t_7=(0,1,1)$\\
$P_8:=(2,4,2):$ $s=3$, $m=7$,  $(0,0,2)$\\
$P_9:=(1,0,2):$ $s=2$, $m=6$,  $\mathbf{ W}=\{(2,4,0),(1,0,2)\}$, $t_9=(1,1,0)$.\\
Then $N:=\{1,x_1,x_2,x_1^2,x_3,x_1^3,x_2x_3,x_3^2,x_1x_2\}$ and so we can easily obtain
$$G=\{x_1^4,x_1^2x_2,x_2^2,x_1x_3,x_2x_3^2,x_3^3\}. $$
The monomials belonging to $G$ are exactly the input for the Axis of Evil algorithm and they are already ordered with respect to our ordering:
\\
starting with $\tau =x_1^4$ we obtain\\
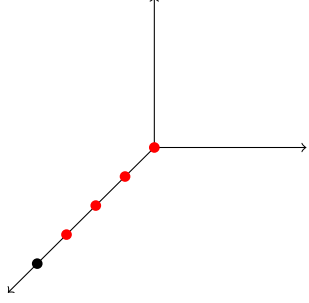
\begin{minipage}[b]{0.5\linewidth}
\centering
\begin{pgfpicture}
\pgfsetarrowsend{to}
\pgfpathmoveto{\pgfpointorigin}
\pgfpathlineto{\pgfpointxyz{0}{0}{5}}
\pgfusepath{stroke}
\pgfpathmoveto{\pgfpointorigin}
\pgfpathlineto{\pgfpointxyz{0}{2}{0}}
\pgfusepath{stroke}
\pgfpathmoveto{\pgfpointorigin}
\pgfpathlineto{\pgfpointxyz{2}{0}{0}}
\pgfusepath{stroke}
\pgfpathcircle{\pgfpointxyz{0}{0}{4}}{2pt}
\pgfusepath{fill}
\color{red}
\pgfpathcircle{\pgfpointxyz{0}{0}{3}}{2pt}
\pgfusepath{fill}
\pgfpathcircle{\pgfpointxyz{0}{0}{2}}{2pt}
\pgfusepath{fill}
\pgfpathcircle{\pgfpointxyz{0}{0}{1}}{2pt}
\pgfusepath{fill}
\pgfpathcircle{\pgfpointxyz{0}{0}{0}}{2pt}
\pgfusepath{fill}
\end{pgfpicture}
\end{minipage}
\hspace{0.5cm}
\begin{minipage}[b]{0.5\linewidth}
$N_1(\tau)=\{1,x_1,x_1^2,x_1^3\}$;\\
$A_1(\tau)=\{(4,0,0),(2,1,4),$\\$(3,0,1),(1,3,4)\}:$  these are the corresponding points via \cemu Correspondence;\\
$B_1(\tau)=\{4,2,3,1\}$\\
$\gamma_{1\tau}=(x_1-4)(x_1-2)(x_1-3)(x_1-1)$: all the linear factors are only depending from $x_1$ are computed in the same time.
\end{minipage}\\
\indent $m=2$:\\
\indent $\zeta_{2\tau}=\gamma_{1 \tau}$\\
\indent  $D_{20}(\tau)=\emptyset$: stop here obtaining,  as  first result, a polynomial   having as leading term an element of $G$ (while the other monomials belong to $N$) and belonging to $I(\mathbf{X})$ since it vanishes in every point of $\mathbf{X}$ (so an element of our minimal Groebner basis).\\
$\tau=x_1^2x_2 $
$N_1(\tau)=\{1,x_1\}$;\\
$A_1(\tau)=\{(2,4,0),(1,0,2)\}$;\\
$B_1(\tau)=\{2,1\}$\\
$\gamma_{1\tau}=(x_1-2)(x_1-1)$
\\
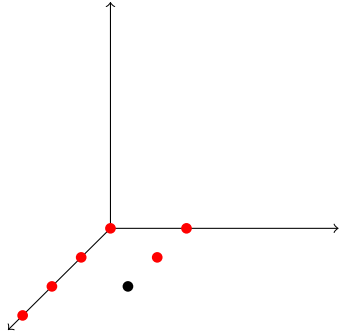
\begin{minipage}[b]{0.5\linewidth}
\centering
\begin{center}
\begin{pgfpicture}
\pgfsetarrowsend{to}
\pgfpathmoveto{\pgfpointorigin}
\pgfpathlineto{\pgfpointxyz{0}{0}{3.5}}
\pgfusepath{stroke}
\pgfpathmoveto{\pgfpointorigin}
\pgfpathlineto{\pgfpointxyz{0}{3}{0}}
\pgfusepath{stroke}
\pgfpathmoveto{\pgfpointorigin}
\pgfpathlineto{\pgfpointxyz{3}{0}{0}}
\pgfusepath{stroke}
\pgfpathcircle{\pgfpointxyz{1}{0}{2}}{2pt}
\pgfusepath{fill}
\color{red}
\pgfpathcircle{\pgfpointxyz{0}{0}{3}}{2pt}
\pgfusepath{fill}
\pgfpathcircle{\pgfpointxyz{0}{0}{2}}{2pt}
\pgfusepath{fill}
\pgfpathcircle{\pgfpointxyz{0}{0}{1}}{2pt}
\pgfusepath{fill}
\pgfpathcircle{\pgfpointxyz{0}{0}{0}}{2pt}
\pgfusepath{fill}
\pgfpathcircle{\pgfpointxyz{1}{0}{0}}{2pt}
\pgfusepath{fill}
\pgfpathcircle{\pgfpointxyz{1}{0}{1}}{2pt}
\pgfusepath{fill}
\end{pgfpicture}
\end{center}
\end{minipage}
\hspace{0.5cm}
\begin{minipage}[b]{0.5\linewidth}
\indent $m=2$:\\
\indent $\zeta_{m\tau}=\gamma_{1\tau}$\\
\indent $D_{20}(\tau)=\{(4,0,0),(3,0,1)\}$\\
We can not stop here, since the obtained polynomial does not vanish at all the points and its head is different from  $\tau \in G$.\\
\indent $N_2(\tau)=\{1,x_1,x_1^2,x_1^3,x_2,x_1x_2\};$
doing so, we find all the monomials of the previous step (we computed their corresponding points) and some new ones.
\end{minipage}
\indent \indent $\delta=1$:\\ 
\indent \indent $A_{21}(\tau)=\{(4,0,0),(3,0,1) \}=D_{20}$ 
\\
The monomials $vx_m^{d_m-\delta}$ are $1,x_1,x_1^2,x_1^3$, corresponding to the points $P_1,P_2,P_4,P_6$.\\
The polynomial already vanishes on $P_2,P_6$, so we consider only the remaining two points.\\
\indent \indent $E_{21}(\tau)=\{1,x_1\}.$ \\
\indent \indent $\gamma_{21\tau}=x_2;$\\
\indent \indent $\xi_{21}=\gamma_{1\tau}\gamma_{21\tau}=(x_1-2)(x_1-1)x_2;$\\
\indent \indent $D_{21}(\tau)=\emptyset$.\\
Remark that $\gamma_{2\tau}$ is actually $\gamma_{21\tau}.$\\
$\tau =x_2^2$\\
$N_1(\tau)=\emptyset$;\\
$A_1(\tau)=\emptyset$;\\
$B_1(\tau)=\emptyset$\\
\indent $m=2$:\\
\indent $D_{20}(\tau)=\mathbf{X}$ \\
\begin{minipage}[b]{0.5\linewidth}
\centering
\begin{pgfpicture}
\pgfsetarrowsend{to}
\pgfpathmoveto{\pgfpointorigin}
\pgfpathlineto{\pgfpointxyz{0}{0}{3.5}}
\pgfusepath{stroke}
\pgfpathmoveto{\pgfpointorigin}
\pgfpathlineto{\pgfpointxyz{0}{2.5}{0}}
\pgfusepath{stroke}
\pgfpathmoveto{\pgfpointorigin}
\pgfpathlineto{\pgfpointxyz{2.5}{0}{0}}
\pgfusepath{stroke}
\pgfpathcircle{\pgfpointxyz{2}{0}{0}}{2pt}
\pgfusepath{fill}
\color{red}
\pgfpathcircle{\pgfpointxyz{0}{0}{3}}{2pt}
\pgfusepath{fill}
\pgfpathcircle{\pgfpointxyz{0}{0}{2}}{2pt}
\pgfusepath{fill}
\pgfpathcircle{\pgfpointxyz{0}{0}{1}}{2pt}
\pgfusepath{fill}
\pgfpathcircle{\pgfpointxyz{0}{0}{0}}{2pt}
\pgfusepath{fill}
\pgfpathcircle{\pgfpointxyz{1}{0}{0}}{2pt}
\pgfusepath{fill}
\pgfpathcircle{\pgfpointxyz{1}{0}{1}}{2pt}
\pgfusepath{fill}
\end{pgfpicture}
\end{minipage}
\hspace{0.5cm}
\begin{minipage}[b]{0.5\linewidth}
\indent $N_2(\tau)=\{1,x_1,x_1^2,x_1^3,x_2,x_1x_2\};$
\indent \indent $\delta=1$:\\ 
\indent \indent $A_{21}(\tau)=\{(2,4,0),(1,0,2) \}$; \\
\indent \indent $E_{21}(\tau)=\{1,x_1\};$\\
\indent \indent $\gamma_{21\tau}=x_2-4x_1+4$\\
\indent \indent $\xi_{21}=\gamma_{1\tau}\gamma_{21\tau}=x_2-4x_1+4;$\\
\indent \indent $D_{21}(\tau)=\{(4,0,0),(2,1,4),(3,0,1),$\\$(2,1,3),(1,3,4)\}$;\\
\indent \indent $\delta=2$:\\
\indent \indent $A_{22}(\tau)=\{(4,0,0),(2,1,4),(3,0,1),$\\$(1,3,4)\}$ 
\\
The terms  $vx_m^{d_m-\delta}$ are  $1,x_1,x_1^2,x_1^3$ and they correspond exactly to $P_1,P_2,P_4,P_6$.
\end{minipage}\\
\indent \indent $E_{22}(\tau)=\{1,x_1,x_1^2,x_1^3\};$\\
\indent \indent $\gamma_{22\tau}=2x_2-x_1^2+7x_1-12;$\\
\indent \indent $\xi_{22}=(x_2-4x_1+4)(2x_2-x_1^2+7x_1-12)$\\
\indent \indent $D_{22}(\tau)=\emptyset$;\\
$\tau=x_1x_3$
\\
\begin{minipage}[b]{0.5\linewidth}
\centering
\begin{pgfpicture}
\pgfsetarrowsend{to}
\pgfpathmoveto{\pgfpointorigin}
\pgfpathlineto{\pgfpointxyz{0}{0}{1.5}}
\pgfusepath{stroke}
\pgfpathmoveto{\pgfpointorigin}
\pgfpathlineto{\pgfpointxyz{0}{1.5}{0}}
\pgfusepath{stroke}
\pgfpathmoveto{\pgfpointorigin}
\pgfpathlineto{\pgfpointxyz{1.5}{0}{0}}
\pgfusepath{stroke}
\pgfpathcircle{\pgfpointxyz{0}{1}{1}}{2pt}
\pgfusepath{fill}
\color{red}
\pgfpathcircle{\pgfpointxyz{0}{1}{0}}{2pt}
\pgfusepath{fill}
\end{pgfpicture}
\end{minipage}
\hspace{0.5cm}
\begin{minipage}[b]{0.5\linewidth}
\centering
$N_1(\tau)=\{1\}$;\\
$A_1(\tau)=\{(2,1,3)\}$;\\
$B_1(\tau)=\{2\}$\\
$\gamma_{1\tau}=(x_1-2)$
\end{minipage}
\\
\indent $m=2$:\\
\indent $N_2(\tau)=\{1\}$.\\
\indent $D_{20}(\tau)=\{(4,0,0),(3,0,1),(1,3,4),(1,0,2)\}$ \\
\indent \indent $\delta=1$:\\
\indent \indent $D_{21}(\tau)=D_{20}(\tau);$\\
\indent $m=3$:\\
\indent $N_3(\tau)=\{1,x_1,x_2,x_1^2,x_3,x_1^3,x_1x_2\}$;\\
\indent $\zeta_{m\tau}=(x_1-2)$;\\
\indent $D_{30}(\tau)=\{(4,0,0),(3,0,1),(1,3,4),(1,0,2)\};$ \\
\\
\begin{minipage}[b]{0.5\linewidth}
\centering
\begin{pgfpicture}
\pgfsetarrowsend{to}
\pgfpathmoveto{\pgfpointorigin}
\pgfpathlineto{\pgfpointxyz{0}{0}{3.5}}
\pgfusepath{stroke}
\pgfpathmoveto{\pgfpointorigin}
\pgfpathlineto{\pgfpointxyz{0}{2}{0}}
\pgfusepath{stroke}
\pgfpathmoveto{\pgfpointorigin}
\pgfpathlineto{\pgfpointxyz{2}{0}{0}}
\pgfusepath{stroke}
\pgfpathcircle{\pgfpointxyz{0}{1}{1}}{2pt}
\pgfusepath{fill}
\color{red}
\pgfpathcircle{\pgfpointxyz{0}{0}{0}}{2pt}
\pgfusepath{fill}
\pgfpathcircle{\pgfpointxyz{0}{0}{1}}{2pt}
\pgfusepath{fill}
\pgfpathcircle{\pgfpointxyz{0}{0}{2}}{2pt}
\pgfusepath{fill}
\pgfpathcircle{\pgfpointxyz{0}{0}{3}}{2pt}
\pgfusepath{fill}
\pgfpathcircle{\pgfpointxyz{1}{0}{0}}{2pt}
\pgfusepath{fill}
\pgfpathcircle{\pgfpointxyz{1}{0}{1}}{2pt}
\pgfusepath{fill}
\pgfpathcircle{\pgfpointxyz{0}{1}{0}}{2pt}
\pgfusepath{fill}
\end{pgfpicture}
\end{minipage}
\hspace{0.5cm}
\begin{minipage}[b]{0.5\linewidth}
\indent \indent $\delta =1$:\\
\indent \indent $A_{31}(\tau)=\{(4,0,0),(3,0,1),(1,3,4),\\(1,0,2)\}$\\
The terms are  $1,x_1,x_1^2,x_1^3,x_2,x_1x_2$, corresponding to $P_1,P_2,P_3,P_4,P_6,\\P_9$,but we can neglect $P_2,P_3$.\end{minipage}\\
\indent \indent $E_{31}(\tau)=\{1,x_1,x_1^2,x_2\};$ \\
\indent \indent $\gamma_{31}(\tau)=6x_3-4x_2+x_1^2-x_1-12;$\\
\indent \indent $\xi_{31}=(x_1-2)(6x_3-4x_2+x_1^2-x_1-12);$\\
\indent \indent $D_{31}(\tau)=\emptyset.$ \\
The desired polynomial is $\gamma_{3\tau}=\gamma_{31}(\tau).$
\\
$\tau=x_2x_3^2$\\
$N_1(\tau)=\emptyset$;\\
$A_1(\tau)=\emptyset$;\\
$B_1(\tau)=\emptyset$\\
\indent $m=2$:\\
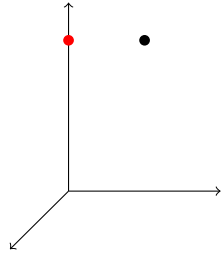
\begin{minipage}[b]{0.5\linewidth}
\centering
\begin{pgfpicture}
\pgfsetarrowsend{to}
\pgfpathmoveto{\pgfpointorigin}
\pgfpathlineto{\pgfpointxyz{0}{0}{2}}
\pgfusepath{stroke}
\pgfpathmoveto{\pgfpointorigin}
\pgfpathlineto{\pgfpointxyz{0}{2.5}{0}}
\pgfusepath{stroke}
\pgfpathmoveto{\pgfpointorigin}
\pgfpathlineto{\pgfpointxyz{2}{0}{0}}
\pgfusepath{stroke}
\pgfpathcircle{\pgfpointxyz{1}{2}{0}}{2pt}
\pgfusepath{fill}
\color{red}
\pgfpathcircle{\pgfpointxyz{0}{2}{0}}{2pt}
\pgfusepath{fill}
\end{pgfpicture}
\end{minipage}
\hspace{0.5cm}
\begin{minipage}[b]{0.5\linewidth}
\indent $N_2(\tau)=\{1\};$\\
\indent $D_{20}(\tau)=\mathbf{X}$;\\
\indent \indent $\delta=1$:\\ 
\indent \indent $A_{21}(\tau)=\{(2,4,2)\};$ \\
\indent \indent $E_{21}(\tau)=\{1\};$\\
\indent \indent $\gamma_{21\tau}=x_2-4$\\
\indent \indent $\xi_{21}=x_2-4;$\\
\indent \indent $D_{21}(\tau)=\{   (4,0,0),(2,1,4),(3,0,1),\\(2,1,3),(1,3,4),(1,0,2)\}$;
\end{minipage}
\\
\indent $m=3$:\\
\indent $\zeta_{3\tau}=x_2-4$\\
\indent $D_{30}(\tau)=D_{21}(\tau)$;\\
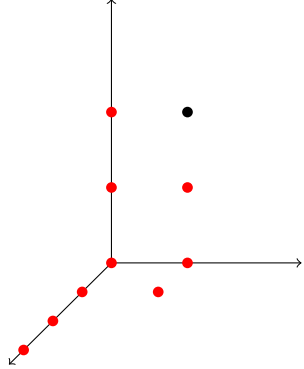
\begin{minipage}[b]{0.5\linewidth}
\centering
\begin{pgfpicture}
\pgfsetarrowsend{to}
\pgfpathmoveto{\pgfpointorigin}
\pgfpathlineto{\pgfpointxyz{0}{0}{3.5}}
\pgfusepath{stroke}
\pgfpathmoveto{\pgfpointorigin}
\pgfpathlineto{\pgfpointxyz{0}{3.5}{0}}
\pgfusepath{stroke}
\pgfpathmoveto{\pgfpointorigin}
\pgfpathlineto{\pgfpointxyz{2.5}{0}{0}}
\pgfusepath{stroke}
\pgfpathcircle{\pgfpointxyz{1}{2}{0}}{2pt}
\pgfusepath{fill}
\color{red}
\pgfpathcircle{\pgfpointxyz{0}{0}{0}}{2pt}
\pgfusepath{fill}
\pgfpathcircle{\pgfpointxyz{0}{0}{1}}{2pt}
\pgfusepath{fill}
\pgfpathcircle{\pgfpointxyz{0}{0}{2}}{2pt}
\pgfusepath{fill}
\pgfpathcircle{\pgfpointxyz{0}{0}{3}}{2pt}
\pgfusepath{fill}
\pgfpathcircle{\pgfpointxyz{1}{0}{0}}{2pt}
\pgfusepath{fill}
\pgfpathcircle{\pgfpointxyz{1}{0}{1}}{2pt}
\pgfusepath{fill}
\pgfpathcircle{\pgfpointxyz{0}{1}{0}}{2pt}
\pgfusepath{fill}
\pgfpathcircle{\pgfpointxyz{0}{2}{0}}{2pt}
\pgfusepath{fill}
\pgfpathcircle{\pgfpointxyz{1}{1}{0}}{2pt}
\pgfusepath{fill}
\end{pgfpicture}
\end{minipage}
\hspace{0.5cm}
\begin{minipage}[b]{0.5\linewidth}
\indent $N_3(\tau)=N(\mathbf{X});$\\
\indent \indent $\delta=1$:	\\
\indent \indent $A_{31}(\tau)=\{(2,1,3)\}$.\\
\indent \indent $E_{31}(\tau)=\{1\};$\\
\indent \indent $\gamma_{21\tau}=x_3-3$\\
\indent \indent $\xi_{31}=(x_2-4)(x_3-3);$\\
\indent \indent $D_{31}(\tau)=\{(4,0,0),(2,1,4),(3,0,1),\\(1,3,4),(1,0,2) \}$;
\end{minipage}\\
\indent \indent $\delta =2$:\\
\indent \indent $A_{32}(\tau)=D_{31}(\tau)$;\\
\indent \indent $E_{32}(\tau)=\{1,x_1,x_1^2,x_1^3,x_2\};$\\
\indent \indent $\gamma_{32\tau}=x_3-4x_2-5x_1^3+41x_1^2-96x_1+48;$\\
\indent \indent $\xi_{32}=(x_2-4)(x_3-3)(x_3-4x_2-5x_1^3+41x_1^2-96x_1+48);$\\
\indent \indent $D_{32}(\tau)=\emptyset;$\\
\indent \indent $\gamma_{3\tau}=(x_3-3)(x_3-4x_2-5x_1^3+41x_1^2-96x_1+48);$ \\
$\tau=x_3^3$\\
$N_1(\tau)=\emptyset$;\\
$A_1(\tau)=\emptyset$;\\
$B_1(\tau)=\emptyset$\\
\indent $m=2$:\\
\indent $D_{20}(\tau)=\mathbf{X}$;\\
\indent $N_2(\tau)=\emptyset;$\\
\indent \indent $\delta=1$:\\ 
\indent \indent $A_{21}(\tau)=\emptyset;$ \\
\indent \indent $D_{21}(\tau)=\mathbf{X}$;\\
\indent $m=3$:\\
\indent $D_{30}=\mathbf{X}$;\\
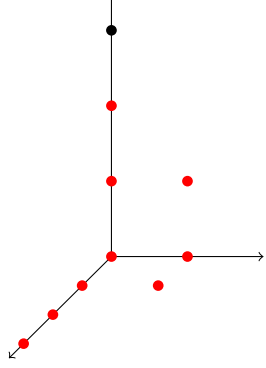
\begin{minipage}[b]{0.5\linewidth}
\centering
\begin{pgfpicture}
\pgfsetarrowsend{to}
\pgfpathmoveto{\pgfpointorigin}
\pgfpathlineto{\pgfpointxyz{0}{0}{3.5}}
\pgfusepath{stroke}
\pgfpathmoveto{\pgfpointorigin}
\pgfpathlineto{\pgfpointxyz{0}{3.5}{0}}
\pgfusepath{stroke}
\pgfpathmoveto{\pgfpointorigin}
\pgfpathlineto{\pgfpointxyz{2}{0}{0}}
\pgfusepath{stroke}
\pgfpathcircle{\pgfpointxyz{0}{3}{0}}{2pt}
\pgfusepath{fill}
\color{red}
\pgfpathcircle{\pgfpointxyz{0}{0}{0}}{2pt}
\pgfusepath{fill}
\pgfpathcircle{\pgfpointxyz{0}{0}{1}}{2pt}
\pgfusepath{fill}
\pgfpathcircle{\pgfpointxyz{0}{0}{2}}{2pt}
\pgfusepath{fill}
\pgfpathcircle{\pgfpointxyz{0}{0}{3}}{2pt}
\pgfusepath{fill}
\pgfpathcircle{\pgfpointxyz{1}{0}{0}}{2pt}
\pgfusepath{fill}
\pgfpathcircle{\pgfpointxyz{1}{0}{1}}{2pt}
\pgfusepath{fill}
\pgfpathcircle{\pgfpointxyz{0}{1}{0}}{2pt}
\pgfusepath{fill}
\pgfpathcircle{\pgfpointxyz{0}{2}{0}}{2pt}
\pgfusepath{fill}
\pgfpathcircle{\pgfpointxyz{1}{1}{0}}{2pt}
\pgfusepath{fill}
\end{pgfpicture}
\end{minipage}
\hspace{0.5cm}
\begin{minipage}[b]{0.5\linewidth}
\indent $N_3(\tau)=N(\mathbf{X});$
\indent \indent $\delta =1$:\\
\indent \indent $A_{31}(\tau)=\{(2,4,2)\};$\\
\indent \indent $E_{31}(\tau)=\{1\};$\\
\indent \indent $\gamma_{31\tau}=x_3-2$;\\
\indent \indent $\xi_{31}=x_3-2$;\\
\indent \indent $D_{31}(\tau)=\{(4,0,0),(2,1,4),(2,4,0),\\(3,0,1),(2,1,3),(1,3,4),(2,4,3)\};$
\end{minipage}\\
\indent \indent $\delta=2$:\\
\indent \indent $A_{32}(\tau)=\{(2,1,3), (2,4,3)\};$\\
\indent \indent $E_{32}(\tau)=\{1,x_2\};$\\
\indent \indent $\gamma_{32\tau}=x_3-3$;\\
\indent \indent $\xi_{32}=(x_3-2)(x_3-3)$;\\
\indent \indent $D_{32}=\{(4,0,0),(2,1,4),(2,4,0),(3,0,1),(1,3,4)\}$;\\
\indent \indent $\delta=3$:\\
\indent \indent $A_{33}(\tau)=D_{32};$\\
\indent \indent $E_{33}(\tau)=\{1,x_1,x_1^2,x_1^3,x_2\}$;\\
\indent \indent $\gamma_{33\tau}=6x_3+8x_2-5x_1^3+35x_1^2-54x_1+24$;\\
\indent \indent $\xi_{33}=(x_3-2)(x_3-3)(6x_3+8x_2-5x_1^3+35x_1^2-54x_1+24);$\\
\indent \indent $D_{33}(\tau)=\emptyset;$\\
The required polynomial is $\gamma_{3\tau}=(x_3-2)(x_3-3)(6x_3+8x_2-5x_1^3+35x_1^2-54x_1+24)$.\\
Then our minimal \gb of the ideal associated to $\mathbf{X}$ with respect to the given order is:
\begin{eqnarray*}
&\mathcal{G}(I(\mathbf{X}))=\Bigl \{x_1^4 - 10x_1^3 + 35x_1^2 - 50x_1 + 24, x_2x_1^2 - 3x_2x_1 + 2x_2,\\ & x_2^2 - 2x_2x_1 - x_2 +  2x^3 - 16x_1^2 + 38x_1 - 24, x_3x - 2x_3 - \frac{2}{3}x_2x_1 + \frac{4}{3}x_2 +\\&+ \frac{1}{6}x^3 - \frac{1}{2}x_1^2 - \frac{5}{3}x_1 + 4, x_3^2x_2 - 4x_3^2 - 7x_3x_2 + 28x_3 + \frac{8}{3}x_2x_1 +\\&+ \frac{20}{3}x_2 - \frac{16}{3}x^3 + 48x^2 - \frac{344}{3}x_1 + 32, x_3^3 - 5x_3^2 + \frac{8}{3}x_3x_2 - \frac{14}{3}x_3 - \frac{16}{9}x_2x_1 \\&-\frac{ 40}{9}x_2+  \frac{73}{9}x_1^3 - \frac{197}{3}x_1^2  +\frac{1358}{9}x_1- 72\Bigr \}, 
\end{eqnarray*}
obtained by our polynomials by the reductions stated in the Axis of Evil Theorem.
\end{example}
Finally, we remark that:
\begin{enumerate}
\item let $\tau_j=x_1^{d_1}\cdots x_n^{d_n}\in G$. The polynomial we are looking for must contain exactly $\sum_{i=1}^n d_i$ factors. It is impossible that the algorithm stops before, so it is impossible that a partial product vanishes on the whole  $\mathbf{X}$. In fact, if so, there would be a polynomial $f \in I$ \st $T(f)\notin (G)$ (we know the minimal basis $G$ \emph{before} starting the Axis of Evil process); 
\item if we otain a factorized polynomial $f$ \st its leading term $T(f)$ belongs to the minimal basis $G$, then $f$ vanishes over all $\mathbf{X}$, because of  \ref{si annulla ovunque}.
\end{enumerate}
\begin{example}\label{Esempio dove non funziona AoE}
Consider the following ideal, given with its primary decomposition:\\
 $J:=(x_1^2,x_2+x_1,x_3) \cap (x_1^2,x_2-x_1,x_3-1)=\\ =(x_1^2,x_1x_2,x_2^2,x_1x_3-\frac{1}{2}x_1-\frac{1}{2}x_2,x_2x_3-\frac{1}{2}x_1-\frac{1}{2}x_2,x_3^2-x_3)\triangleleft \CC[x_1,x_2,x_3]$.\\ 
Call its generators $f_1,...,f_6$, considering them in the correct order.\\
It is $0$-dimensional because $x_1^2,x_2^2,x_3^2\in In(J)$ (see \cite{SPES2}), but it is not radical: its radical is  $\sqrt{J}=(x_2,x_3^2-x_3, x_1).$\\
For such an ideal the Axis of Evil does not hold.\\
Consider the polynomial $f_4=x_1x_3-\frac{1}{2}x_1-\frac{1}{2}x_2.$\\
By the Axis of Evil theorem (\ref{ASSE DEL MALE}), its factorization should be of the form:
$$(x_1+...)(x_3+...) $$
and we should have
$$x_1x_3-\frac{1}{2}x_1-\frac{1}{2}x_2+Px_1^2+Qx_1x_2+Rx_2^2,\; P,Q,R \in \CC[x_1,x_2,x_3],$$
since we can only reduce deleting the multiples of $x_1^2,x_1x_2,x_2^2$, in order to obtain $f_4$. In order to have the correct product we must have $-\frac{1}{2}x_2$ in it.
We can not obtain it through reductions, so the only chance is 
that we have a product of the form
$$k*hx_2, $$
with $h,k$ constants \st $hk=-\frac{1}{2}$, in particular both different from $0$.\\
A priori, we can have two possibilities:
\begin{itemize}
\item $(x_1+k+...)(x_3+hx_2+...)$;
\item $(x_1+hx_2+...)(x_3+k+...)$.
\end{itemize}
The second one is impossible: the polynomial having $x_1$  as head can not contain variables greater than $x_1$, so we consider only:
$$(x_1+k+...)(x_3+hx_2+...).$$
We will then obtain
$$x_1x_3+hx_1x_2+kx_3-\frac{1}{2}x_2+...$$
We can delete the term $x_1x_2$ but it remains $kx_3$ which can not be reduced.
\end{example}
\section{Corollaries.}
We enumerate here some famous theorems which can be easily proved as corollaries of the Axis of Evil Theorem. For more details see, for example, \cite{SPES2}.\\
Here we provide the general statements of these results, but clearly they can only be deduced under the hypothesis of the Axis of Evil theorem
\\
The first one is \emph{Lazard Structural Theorem}, which describes the structure of a minimal lexicographical \gb of an  $I \triangleleft k[x_1,x_2].$
\\
The original proof considers $P=k[x_1,x_2]=k[x_1][x_2]$ and it is based on the fact that $k[x_1]$ is a Principal Ideal Domain (PID).\\ 
Norton-S\u al\u agean \cite{NS} reformulated it using, more generally, $R[x]$ with $R$ PIR.\\
We briefly recall the following
\begin{definition}
The \emph{content} $r_{f}\in R$, with $R$ PIR, of a polynomial $f(x) \in R[x]$ is the GCD of its coefficients. A polynomial   $f(x) \in R[x]$ is called \emph{primitive} if $r_f=1$.\\
The \emph{primitive part} of $f(x) \in R[x]$ is the polynomial $p_0(x) \in R[x]$ \st $$f(x)=r_fp_0(x).$$
\end{definition}
Let $R$ be a PIR, $P:=R[x]$. Let $I \triangleleft P$ e $F:=\{f_0,...,f_s\}$ a \emph{minimal} \gb of $I$ ordered in such a way that, called $d(i):=deg(f_i),$  $\forall i$, $0 \leq i \leq s$
$$d(0)\leq ...\leq d(s). $$
Define then $c_i=lc(f_i),\, r_i \in R \setminus \{0\}$ e $p_i \in P$ the leading coefficient, the content and the primitive part of $f_i$, for all $1\leq i \leq n$.
\begin{theorem}[Lazard]\label{SPES 33.1.5}
If, moreover, $R$ is a PID, then:
\begin{itemize}
\item $f_0=PG_1 \cdots G_{s+1};$
\item $f_j=PH_jG_{j+1}\cdots G_{s+1},\, 1\leq j \leq s.$
\end{itemize}
where
\begin{enumerate}
\item $d(1)<...<d(s);$
\item $G_i \in R,\; 1 \leq i \leq s+1$ is \st $c_{i-1}=G_ic_i$ 
\item $P=p_0$ (the primitive part of $f_0\in R[x]$);
\item $H_i \in R[x]$ is a monic polynomial of degree  $d(i)$ in $x$, for all $i$;
\item for all $i$ we have $H_{i+1}\in (G_1 \cdots G_i,H_1G_2\cdots G_i,...,H_{i-1}G_i,H_i);$
\item $r_i=G_{i+1}\cdots G_s$ 
\end{enumerate}
\end{theorem}
\begin{theorem}[Norton-S\u al\u agean]\label{NS}
With the previous notation, each $$p_i \in (f_j, j<i):r_i.$$
\end{theorem}
In fact, we have $r_i=\prod_{m=1}^{n-1} \prod_{\delta=1}^{d_m}\gamma_{m\delta t_i}$ and $p_i= \prod_{\delta=1}^{d_n}\gamma_{n\delta t_i}$.
\\
The second well-known result which can be straightforwardly derived from the Axis of Evil Theorem is the well known \emph{Elimination Theorem} (see \cite{B} for details)
\begin{theorem}[\cite{Tr}]
Let $I\triangleleft k[x_1,...,x_n]$ an ideal, take the lexicographical ordering induced by $x_1<...<x_n$ and call $I_j$ the $j$-th elimination ideal $I_j=I\cap k[x_1,...,x_j]$. Let $\mathcal{G}$ be a \gb of $I$, then $\mathcal{G}_j=\mathcal{G}\cap k[x_1,...,x_j]$ is a \gb of $I_j$.
\end{theorem}
The following result, \emph{Kalkbrener theorem} (\cite{K2}, \cite{SPES2}), is another consequence of the Axis of Evil Theorem and it is a stronger characterization of the lexicographical ordering.\\
For each subset $L\subset k[x_1,...,x_n]$, $i=1,...,n$, $\forall \delta \in \NN$ set 
$$L_{i\delta}=\{ p \in L,\,\vert p \in k[x_1,...,x_i], \,deg_i(p)\leq \delta\} $$
and
$$Lp_{i, \delta}=\{  Lp(p),\,p \in L_{i, \delta}\}.$$
\begin{theorem}[Kalkbrenner]
With the previous notations, considered an ideal $I \triangleleft k[x_1,...,x_n]$ and a \gb $\mathcal{G}$ of it, these forms are equivalent:
\begin{itemize}
\item $\mathcal{G}$ is a \gb of $I$ w.r.t, the lexicographical order $<$ induced by $x_1<...<x_n$;
\item $Lp_{i,\delta}(\mathcal{G})$ is a \gb of $Lp_{i, \delta}(I)$,  $i=1,...,n$, $\forall \delta \in \NN$.
\end{itemize}
\end{theorem}
Let us now mention \emph{Gianni-Kalkbrener theorem}, whose situation is a bit more complicated (see \cite{K}, \cite{G}, \cite{SPES2}).\\
\begin{theorem}[Gianni-Kalkbrener]
Let $I \triangleleft k[x_1,...,x_n]$ an ideal and $\mathcal{G}$ w.r.t the lexicographical order $<$ induced by $x_1<....<x_n$. As before we define also $\mathcal{G}_d=\mathcal{G}\cap k[x_1,...,x_d]$.\\
Consider $\alpha=(b_1,...,b_d)\in V(I_d)$ and define the projection map
$$\Phi_{\alpha}:k[x_1,...,x_n]\rightarrow k[x_{d+1},...,x_n] $$
$$f(x_1,...,x_n) \mapsto f(b_1,...,b_d,x_{d+1},...,x_n).$$
Let $\sigma$ be the minimal value such that $\Phi_{\alpha}(Lp(g_{\sigma}))\neq 0$ and $j,\delta$ the values \st
$$g_{\sigma}=Lp(g_{\sigma})x_j^{\delta +1}+... \in k[x_1,...,x_j] \setminus k[x_1,...,x_{j-1}].$$
Then
\begin{enumerate}
\item $j=\delta +1$
\item $\forall g \in \mathcal{G}_d$, $\Phi_{\alpha}(g)=0$;
\item $\forall g \in \mathcal{G}_{d+\delta}$, $\Phi_{\alpha}(g)=0$;
\item $\Phi_{\alpha}(g_{\sigma})=gcd(\Phi_{\alpha}(g),\, g \in G_{d+1})\in k[x_{d+1}]$;
\item $\forall b \in k$, $(b_1,...,b_2,b) \in V(I_{d+1}) \Leftrightarrow \Phi_{\alpha}(g_{\sigma})(b)=0.$
\end{enumerate}
\end{theorem}
Clearly $(1-3)$ are essentially a corollary of theorem \ref{NS}; on the other side, $(4-5)$ apparently cannot be deduced from the Axis of Evil Theorem.
\section{Acknowledgement.}
I wish to thank M. G. Marinari for her help, ideas and suggestions while studying this subject.
\end{document}